\documentclass[14pt]{extarticle}
\usepackage[margin=1in]{geometry}
\usepackage{dsfont}
\usepackage{amsfonts}
\newtheorem{theorem}{Theorem}[section]
\newtheorem{pro}{Proposition}[section]
\newtheorem{lemma}{Lemma}[section]
\newtheorem{definition}{Definition}[section]
\newtheorem{remark}{Remark}[section]
\newtheorem{cor}{Corollary}[section]

\newcommand{\proof}[1]{\noindent{\it\bf Proof:#1\ }}
\newcommand{\QED}{\hfill$\Box$\medskip}

\begin{document}

\title{On  Stability of  Nodal $L_k^p$-Maps I}
\author{Gang Liu\\Department of Math, UCLA  }
\date{}
\maketitle

\section{Introduction}
In this paper and its sequels, we  introduce the  concept   of  weak stability for  nodal $L_k^p$-maps as  a natural generalization of  stability for    stable $J$-holomorphic maps in Gromov-Witten theory introduced by Kontsevich in \cite{2}. We 
 then give a complete characterization of  the weakly stable  nodal $L_k^p$-maps in term of  their isotropy   groups.
  Among  weakly stable $L_k^p$-maps, the stable ones are those whose isotropies are finite.  As a consequence, we prove   that  the space of unparametrized
  nodal $L_k^p$-maps modeled on a fixed tree is always  Hausdorff  without any stability
  assumption.
   
   We will only deal with the genus zero case. This is justified since for a  nodal $L_k^p$-map, the  part of its  reparametrization group with positive dimension consists of the reparametrizations of  its  unstable genus zero components. In this paper, we only consider the  nodal $L_k^p$-maps modeled on a fixed tree $T$. The general cases
allowing the changes of the topological types of the domains   and   targets will be treated in  the sequels of this paper. We now describe the main results of this paper using the terminologies defined in next section.

Fix a (minimal) label  $L$ of a given tree $T$ such that the labeled
tree ${\hat T}=(T, L)$ is stable. Let ${\cal M}_{\hat T}$ be the moduli space of  genus zero stable  curves with $n$ marked points modeled on ${\hat  T}$ and ${\cal U}_{\hat T}\rightarrow {\cal M}_{\hat T}$ be the
universal curve over
${\cal M}_{\hat T}$. Here $n$ is the minimal number of marked points added to a genus zero nodal surface modeled on $T$ to make it stable. Fix a
Riemannian manifold $M$ of dimension at least $2.$
Let ${\tilde {\cal B}}^{\hat T}={\tilde {\cal B}}^{\hat T}_{k, p}$  be the set of $L_k^p$-maps $f:\Sigma\rightarrow M$ with $\Sigma$ being one of the fibers of the universal  family  ${\cal U}_{\hat T}\rightarrow {\cal M}_{\hat T}$. One can show that ${\tilde {\cal B}}^{\hat T}_{k, p}$ is a Banach  manifold of class $[m_0]$ (see Section 2).  Here $m_0=k-2/p$ is the Sobolev differentiability of $f$.  Throughout this paper, we will alway assume that
 $[m_0]\geq 1$  so that each component of $f$ is at least of class $C^1.$ Roughly speaking, the space $ {\tilde {\cal  B}}^{\hat T}$ can be thought as the space of parametrized nodal $L_k^p$-maps modeled on $T$. Let $  {\cal  B}^{\hat T}$ be the space of equivalence classes of nodal $L_k^p$-maps. It can be obtained from $ {\tilde {\cal  B}}^{\hat T}$ as the  orbit space under the actions of the compatible system of reparametrization groups $\{G_f|\, f\in  {\tilde {\cal  B}}^{\hat T}\}$. The normal subgroup  of $G_f$   preserving the components of $f$  is independent  of $f$ in the sense that they can be identified with each other
 canonically up to an overall conjugation. Denote the resulting group by $G_T$. Note that the quotient group $G_{f/T}=G_f/G_T$ is  a finite group that exchanges  the components of $f$. Thus $  {\cal  B}^{\hat T}$ can be obtained by first forming the global quotient $  {\tilde {\cal  B}}^{\hat T}/G_{T}$, then   quotient out a further locally finite equivalence relation by the actions of $G_{f/T}.$
 Note that the action of $G_T$ is continuous (see Sec. 2).

 \begin{theorem}
 	The space  $   {\tilde {\cal  B}}^{\hat T}$ is $G_T$-Hausdorff in the sense that for any two different $G_T$-orbits $G_Tf_1$ and $G_Tf_2$, there exit $G_T$-neighborhoods $G_TU_1$ and $G_TU_2$ such that $G_TU_1 \cap G_TU_2=\varphi$. Therefore, the global quotient $  {\tilde {\cal  B}}^{\hat T}/G_{T}$ is Hausdroff.
 \end{theorem}

 The proof of the above theorem implies  the following

 \begin{pro}
The space $  {\cal  B}^{\hat T}$ of the unparametrized nodal $L_k^p$-maps  is always  Hausdroff.

\end{pro}

 Note that in general a nodal map  $f$ may have components, such as trivial unstable component, so that the istropy group ${\Gamma_f}$ is
 not compact. In this case, the action of $G_T$ on $   {\tilde {\cal  B}}^{\hat T}$ is certainly not proper so that $f$ is not stable in any reasonable sense. Yet, the above result shows that when the topological type of the domains  is fixed given by $T$, the Hausdorffness of $  {\tilde {\cal  B}}^{\hat T}/G_{T}$  and  $  {\cal  B}^{\hat T}$ still holds  without  requiring any stability conditions. This seems to contradict to our experience in Gromov-Witten theory. Indeed, when
 the topological type of the domains  is allowed to change, the corresponding space ${\cal B} $ of unparemetrized nodal $L_k^p$-maps is not Hausdroff anymore. But this non-Hausdorffness occurs in a rather definite manner mainly caused by the appearance of the extra trivial bubbles obtained by a non-convergence sequence of reparametrizations. Once such  degenerations are prohibited, the-Hausdorffness will be restored even allowing the change of the topological types of the domains. In other words, for any subspace of  ${\cal B} $ which do not contain a sequence convergent to a nodal map with extra trivial bubbles, the Hausdorffness still holds.  

 Next we define the weak stability.

 \begin{definition}
  A nodal $L_k^p$-map is said to be weakly stable if none of its unstable components is a trivial map.	
 \end{definition}

 Let  ${\tilde {\cal  B}}^{ws}_{\hat T}$ be the subspace  of $ {\tilde {\cal  B}}^{\hat T} $  consisting of weakly stable  nodal $L_k^p$-maps. As before,   ${\tilde {\cal  B}}^{ws}_{\hat T}$ can be thought as the space of parametrized  weakly stable nodal maps modeled on $T$. Let  ${\cal  B}^{ws}_{\hat T}$ be the corresponding space of
 unparametrized  weakly stable nodal maps.

\medskip
\noindent
 {\bf  Note:}   Applying  the above definition to   $J$-holomorphic  nodal maps, we get
   one of the standard definitions for  stable $J$-holomorphic maps
 in Gromov-Witten theory \cite{2}. 
  In this sense the weak stability for $L_k^p$-maps here is a natural generalization
of the stability for  $J$-holomorphic maps.  

Recall in the case that $X$ is a locally compact topological space such as a finite dimensional manifold and $G$ is a Lie group or locally compact topological group, a (continuous) group action $\Phi:G\times X\rightarrow X$ is said to be proper if the map  $\Phi\times{\it id}_{X}:G\times X\rightarrow X\times X$  is a proper map:  the inverse image $(\Phi\times{\it id}_{X})^{-1}(K) $ of any compact subset $K\subset X\times X$ is compact.

 In our infinite dimensional case, $X=:   {\tilde {\cal  B}}^{ws}_{\hat T}$
is not locally compact, the above definition is too weak to be useful.

We introduce the following  stronger definition.

\begin{definition}
A group action $\Phi:G\times X\rightarrow X$ is said to be proper if  for any compact subset $K\subset X\times X$, there is a neighborhood $U$ of $K$ in  $ X\times X$ such that the image  $\pi_G ((\Phi\times{\it id}_{X})^{-1}(U))$ of the projection to $G$ of the inverse image $(\Phi\times{\it id}_{X})^{-1}(U)$  is pre-compact in $G$.

\end{definition}


In Sec. 3, we will show that if $G$ and  $X$ are locally compact,
the definition here is   equivalent to the usual definition.

 \begin{theorem}
 	The action of $G_T$ on  $   {\tilde {\cal  B}}^{ws}_{\hat T}$ is
 	proper.

 \end{theorem}



 \begin{cor}

 For any weakly stable nodal map $f$,  the  isotropy $\Gamma_f$ of the $G_T$-action or $G_f$-action  is always compact.
 \end{cor}


 {\bf  Note:}
 It follows from the continuity of the $G_T$-action that   $\Gamma_f$  is closed in  $G_T$ so that it is a  compact  Lie subgroup of $G_T$ if $f$ is weakly stable.

 \begin{cor}
 	
 	A  nodal $L_k^p$ map $f$ is weakly stable if and only if the  isotropy $\Gamma_f$ is compact.

 \end{cor}

\medskip
This gives a complete characterization of weak stability of  nodal $L_k^p$-maps in term of their  isotropy groups.

\medskip
 \begin{definition}
 	A weakly stable nodal $L_k^p$-map $f$ is said to be  stable if
 	it has no infinitesimal automorphisms.  In other words, the dimension
 	of the Lie algebra  $L(\Gamma_f)$ is equal to zero.
 \end{definition}

 Since $\Gamma_f$ is compact for a  weakly stable nodal map, the above condition is equivalent to the following definition.

\begin{definition}
	A weakly stable nodal $L_k^p$-map $f$ is said to be  stable if
its isotropy group 	$\Gamma_f$ is a finite group.
\end{definition}

Using the above characterization  for a weakly stable nodal $L_k^p$-map, we get
 the following corollary.

\begin{cor}
	
	A   nodal  $L_k^p$-map $f$ is stable if and only if the  isotropy $\Gamma_f$ is finite.
	
\end{cor}

Recall that the corresponding characterization  for stable $J$-holomorphic maps\cite{2} is the following well-known proposition in GW theory.

\begin{pro}
	
	A $J$-holomorphic  nodal  map $f$ is stable if and only if the  isotropy $\Gamma_f$ is finite.
	
\end{pro}

 Thus stable $L_k^p$-maps and stable $J$-holomorphic    maps have exactly the same
 characterizations in term of their isotropies. One the other hand, we have seen that
 in term of the definitions,  weakly stable $L_k^p$-maps are the natural generalizations of   stable $J$-holomorphic   maps.


  \medskip
  \medskip

 Let
  $   {\tilde {\cal  B}}^{s}_{\hat T}$ be the  set  of stable nodal $L_k^p$-maps and $    {\cal  B}^{s}_{\hat T}$ be the corresponding
  quotient space of unparametrized stable nodal maps.  Next theorem follows from the definition.

  \begin{theorem}
  	The actions of $G_T$ and $G_f$  on  $   {\tilde {\cal  B}}^{s}_{\hat T}$  are
  proper with finite istropies so that $   {\cal  B}^{s}_{\hat T}$ is
   a topological Banach orbifold.    In particular, the subspace
      of  $   {\cal  B}^{s}_{\hat T}$ consisting of stable nodal
      maps with the  trivial isotropy is a topological Banach manifold.
  	
  \end{theorem}

  To get a better  understanding of the stability, we divide the unstable  components of  weakly stable maps into the following two classes.
  Let $f_v:\Sigma_v\simeq {\bf C P }^1\rightarrow M$ be such a component.
  Then $f_v$ is said to be $2$-dimensional if there is a point $x_0$ in
  $\Sigma_v$ such that the rank of $(df_v)_{x_0}$ is equal to two; otherwise it is one-dimensional.

 The  following proposition is proved in Sec. 4. 

\begin{pro}
	If all unstable components of a weakly stable map are
	$2$-dimensional,  it is stable.
\end{pro}

\begin{remark}
(1) In above proposition, the condition that all unstable components of $f$
are
$2$-dimensional already implies that all such  components
are nontrivial, hence $f$ is  weakly stable by definition.

(2)	As a comparison,  in \cite{3} the notion of stability for  nodal $L_k^{p}$-maps $f$ with $\delta$-exponential decay  in  a symplectic manifold is defined by requiring
	that the symplectic area of any unstable  component of $f$ is positive. 
	Note that the positivity of the symplectic area of any unstable  component implies that it is  $2$-dimensional in the definition above, hence $f$ is stable.
\end{remark}

In last section of this paper, we will discuss the case of a weakly stable map
$f$ with $1-$dimensional unstable components and prove the following proposition (see the definition in the last section).

\begin{pro}
	If none of the  unstable  components of a weakly stable map is standard $S^1$-invariant, then it is stable.
	
\end{pro}

Since the component that is  standard $S^1$-invariant can be explicit identified (up to a conjugation), this proposition
 means that  other than   this exceptional case all  weakly stable maps
with $1-$dimensional unstable components are still stable.

Of course, if an unstable  component of a weakly stable map $f$ is standard $S^1$-invariant, its isotropy group is $S^1$ so that $f$ is not stable but only weakly stable.

 More generally,  we have the  next proposition  on the continuous part of the isotropy
  group of  a nodal map.

  \begin{pro}
  If the identity component  $\Gamma_f^0$ of the isotropy group of the nontrivial components of a	 nodal $L_k^p$-map  $f$ is nontrivial, it is a torus $T^n=(S^1)^n$. In particular, $\Gamma_f^0\simeq T^n$ for a weakly stable but not stable  $L_k^p$-map $f.$ This can only happen if $f$ is  standard $S^1$-invariant (up to a conjugation).
  \end{pro}

  This paper is organized as follows.

  Section 2 provides the  preliminaries used to define the space ${\tilde {\cal B}}_{\hat T}$ of $L_k^p$-maps modeled on a labeled tree ${\hat T}.$

  Section 3 proves the first two   theorems  here as well as  some related results.
  The results here generalize the corresponding ones in \cite{6} for $L_k^p$-maps with  domain $S^2$.

   Section 4 discusses the isotropy groups of the  weakly stable $L_k^p$-maps.
 
This paper is a revised version of the paper with  the same title \cite{7}. The main ideas were already in \cite{6}.

   \section{Preliminaries}

In the first two  subsections, we recall the  construction of    the moduli space ${\cal M}_{{\hat T}}$ of stable curves modeled on a fixed $n$-labeled tree ${\hat T}$. The  results are well-known. We refer the readers to \cite{9} and the reference therein for details.

 \subsection{Nodal curves, stable curves  modeled on ${\hat T}$}

A nodal curve/surface  ${\bf \Sigma}$ is an  union of finitely many its components, with each component being  a compact Riemann surface  (as $1$-dimensional complex manifold), joining together at their double points.
Throughout this paper, we always assume that ${\bf \Sigma}$ is connected.

Given  a nodal curve ${\bf \Sigma}$,  the  associated graph $T_{\bf \Sigma}$ encodes its combinatoric information of its components and is defined  as following: each component $\Sigma_v$ corresponds to a vertex $v$ of $T_{\bf \Sigma}$ and a double point $d_{uv}=d_{vu}$ lying in the two components $\Sigma_u$ and $\Sigma_v$ corresponds to an edge denoted by $[uv]=[vu] $.

In the following we will
abuse notation and denote the set of vertices by $T_{\bf \Sigma}$ as well.
Denote the edge set of a graph $T$ by $E_T$.

A genus zero  nodal curve ${\bf \Sigma}$ is a nodal curve such that its
``arithmetic genus" is equal to zero; equivalently  each component of
${\bf \Sigma}$ is a smooth genus zero curve and the  associated graph $T_{\bf \Sigma}$ is a tree. 

From now on all nodal curves or $n$-marked nodal curves defined below  will always have genus zero.

In this paper, we  fix a tree $T$  and consider  nodal curves ${\bf \Sigma}$ with $T_{\bf \Sigma}=T$. Each such curve is said to be modeled on $T$.
Thus in this notation, a  nodal curve ${\bf \Sigma}$
 modeled on $T$ is the pair ${\bf \Sigma}=(\Sigma, {\bf d})$ with underlying surface $\Sigma=\coprod_{v\in T}\Sigma_v/\{{d_{uv}=d_{vu}, [uv]}\in E_T\}.$ Here $\Sigma_v\simeq S^2$ is the component associated to   the vertex $v\in T$, and ${\bf d}=\cup_{v\in T} {\bf d}_v$ with   ${\bf d}_v=\{d_{uv}, [uv]\in E_T\}$  being the set of double points on $\Sigma_v$.

An $n$-marked  nodal curve is a triple ${\bf \Sigma}=(\Sigma, {\bf d}, {\bf x})$,
 where  $(\Sigma, {\bf d})$ is a  nodal curve and  $ {\bf x}=\cup_{v\in T}{\bf x}_v$, where ${\bf x}_v=\{x_{vi}\}$ is the set of (distinct) marked points on $\Sigma_v$. We require that  ${\bf x}\cap  {\bf d}=\varphi.$
  If  the number of the marked points ${\bf x}$ is $n$, the combinatoric information on how these marked points lying the components of  $\Sigma$  can be encoded by a map $L: {\bf n}=\{1, 2, \cdots, n\}\rightarrow
 T_{(\Sigma, {\bf d})}$, which at the same time  gives the order of the
 $n$-marked points.

 An $n$-labeled tree ${\hat T}=(T, L)$  consists of   a tree $T$ and a  label $L$
  as above. 
  
  From now on, we will fix a $n$-labeled tree ${\hat T}$ and consider all $n$-marked genus zero nodal curves  ${\bf \Sigma}$ (modeled on ${\hat T}$) with ${ \hat T}_{\bf \Sigma}={\hat T}$.

Let ${\bf p}_v= {\bf d}_v\cup  {\bf x}_v $ be set of special points on $\Sigma_v$.  A $n$-marked    curve ${\bf \Sigma}$ is said to be stable if on each component  $\Sigma_v$,  the number of special points, $\#({\bf p}_v)\geq 3.$

A component  of  a  nodal curve or a $n$-marked    curve is said to be  unstable if it has  at most   two double points.  A unstable component ${\Sigma}_v, v\in T$ with only one (or without ) double point  will be call a top bubble (on the bubble tree $T$ with $v$ being a tip of $T$).

  Given a genus zero nodal curve $(\Sigma, {\bf d})$, a stable curve $(\Sigma, {\bf d}, {\bf x})$, called (one of) the  minimal stabilization(s),   can be obtained by adding   minimal number of marked points ${\bf x}$ to the unstable components, which is well defined up to  automorphisms  of $\Sigma$ and a choice of the labeling of the marked points ${\bf x}$. For the purpose  of defining the space ${\tilde {\cal  B}}^{\hat T}$   parametrized
  nodal maps below, 
 it is sufficient to  only consider stable curves with minimal number of marked points.

  A map $\phi: \Sigma_1\rightarrow \Sigma_2$ is said to be an equivalence between two nodal curves  ${\bf \Sigma}_1=(\Sigma_1, {\bf d}_2)$ and ${\bf \Sigma}_2=(\Sigma_2, {\bf d}_2)$ (modeled on $T_1$  and $T_2$ respectively)
    if it is a homeomorphism such that
    for each $v\in T_1$, $\phi_v:(\Sigma_1)_v\rightarrow (\Sigma_2)_{\phi (v)}$ is biholomorphic. In the case that $\Sigma_1= \Sigma_2=\Sigma,$ the equivalence $\phi$ above is called a automorphism of ${\bf \Sigma}.$

   If in addition, ${\bf \Sigma}_1=(\Sigma_1, {\bf d}_2, {\bf x}_1)$ and ${\bf \Sigma}_2=(\Sigma_2, {\bf d}_2, {\bf x}_2)$ (modeled on ${\hat T}_1$  and ${\hat T}_2$ respectively) are $n$-marked nodal curves, then
  	 a map $\phi: \Sigma_1\rightarrow \Sigma_2$ is said to be an equivalence if it is an equivalence for  the underlying nodal curves that  commutes with the labeling $L_1$ and $L_2$.

  \begin{lemma}
  Let $(\Sigma_i,{\bf d}_i, {\bf x}_i)$ be a minimal stabilization of a nodal curve  $(\Sigma_i,{\bf d}_i), $ $ i=1, 2 $ and
 $\phi: \Sigma_1\rightarrow \Sigma_2$ is an equivalence  between two nodal curves. Then $\phi$ is determined by its effect on  ${\bf x}^t_1$, where  ${\bf x}^t_1\subset  {\bf x}_1$  is the set of the  marked points that are lying on the top bubbles.	
  \end{lemma}
 \proof

 The proof is a simple induction on the number of top bubbles. We leave it the  readers.

 \QED

  Applying this to an automorphism of  a stable   curve, we get the following well-known corollary.

 \begin{cor}
 Let  $\phi:  \Sigma\rightarrow \Sigma$ be  an  automorphism of  a stable   curve  ${\bf \Sigma}=(\Sigma, {\bf d}, {\bf x})$. Then $\phi=id.$	
\end{cor}

  Given a nodal curve $(\Sigma, {\bf  d})$ model on $T$, denote the group of its automorphisms by   $G_{\Sigma}$.   Let  $G_{T}$ be the continuous part  of
  $G_{\Sigma}$, which  is the subgroup consisting of automorphisms that map each  component of $\Sigma$ to itself.  Thus each element $\phi \in G_{T}$ acts trivially on a
   stable component, and on each unstable component $\Sigma_v$, $\phi=\phi_v$ is an automorphism of $\Sigma_v$ that fixes the double point(s).

  Clearly
  $G_{\Sigma}$ depends on $\Sigma$ in general. However, $G_{T}$ only depends on $T$ (up to a conjugation).

  To see this, fix a minimal stabilization  $(\Sigma, {\bf  d}, {\bf x})$ modeled on a $n$ labeled tree ${\hat T}$. For each unstable component $(\Sigma_v, {\bf  d}_v, {\bf x}_v)$,  its  three     double/marked  points are ''canonically" ordered by the process in next subsection. Hence  we get a
   well-defined identification  $\Sigma_v\simeq {\bf CP}^1_v={\bf CP}^1$ by sending above three points to the points $(0, 1, \infty)$ in ${\bf CP}^1$.
 By conjugating by this identification,
 the reparametrization group $G_v=G_{\Sigma_v}$ of $\Sigma_v$  can be identified with one of the
  subgroups
   $G_i, i=0, 1,2$ of $G_0={\bf PSL}(2, {\bf C})$ consisting of the biholomorphic maps of ${\bf CP}^1_v$ that fix the first $i$ points of the three standard points.

  Then  $G_T$ is isomorphic to  the product $\prod_v G_v$ over   unstable components $\Sigma_v$.  Of course the isomorphism is only canonical up to a conjugation due to the choice of the above identification.

   We now introduce a few  finite subgroups of $G_{\Sigma}$. Let  $(\Sigma, {\bf  d}, [{\bf  x}])$ be the   stabilization of a nodal curve
   $(\Sigma, {\bf  d})$ by adding minimal number of unordered marked points $[{\bf x}]$, and $(\Sigma, {\bf  d}, [{\bf  x}]^0)$ be the   stabilization of 
   $(\Sigma, {\bf  d})$ by give a  partially ordering to the unordered marked points $[{\bf x}]$.  Here the partial ordering of $[{\bf x}]^0$ is  defined as following:  the two marked points on each top bubble are ordered but no other ordering among the marked points $[{\bf x}]$.

   This minimal stabilization is unique 
   up to an equivalence.
   Here the equivalence of two marked nodal curves with unordered or partially ordered  markings is defined
  in obvious way similar to the ordered case.  Denote  the group of  automorphisms (=self equivalence) of $(\Sigma, [{\bf  x}])$ and $(\Sigma, [{\bf  x}]^0)$ by  $G_{\Sigma, [{\bf  x}]}$ and $G^0_{\Sigma, [{\bf  x}]}$ respectively. Unlike the  ordered case, $G_{\Sigma, [{\bf  x}]}$ and $G^0_{\Sigma, [{\bf  x}]}$  are  not trivial but  a  finite group. By the above lemma, each element of $G_{\Sigma, [{\bf  x}]}$ is  given by a permutation of the set of marked points, and 
  $G^0_{\Sigma, [{\bf  x}]}$  is  the subgroup of  $G_{\Sigma, [{\bf  x}]}$ consisting the automorphisms that preserve the order of the two marked points of each top bubble. In particular, any non-trivial element in  $G^0_{\Sigma, [{\bf  x}]}$
   switches the components of $\Sigma.$

  Since each element of $G_T$ preserve all the components, the intersection
 $G_T\cap  G^0_{\Sigma, [{\bf  x}]}$ is trivial. Finally  the finite  group
 $G_{\Sigma}/G_T$ can be considered as a subgroup of the group  that  permutes the components of $\Sigma.$

  \begin{lemma}
  	The group $G_T$ is a  normal subgroup of $G_{\Sigma}.$
  	 The exact sequence $1\rightarrow G_T\rightarrow 	 G_{\Sigma}\rightarrow 	G_{\Sigma}/G_T\rightarrow 1$  splits by $G^0_{\Sigma, [{\bf  x}]}$.

  	 Moreover $G_T$ and $G^0_{\Sigma, [{\bf  x}]}$ generate 	$G_{\Sigma}$ and the induced homomorphism $G^0_{\Sigma, [{\bf  x}]}\rightarrow G_{\Sigma}/G_T$ is an isomorphism.
  Consequently  $G_{\Sigma}$ is  the semi direct product of $G_T$ and $G^0_{\Sigma, [{\bf  x}]}$.

  \end{lemma}

  \proof


  Since any element $\phi$ of $G_{\Sigma}$ is determined by its effect on the marked points on the top bubbles, by composing an element $\psi$ in $G_T$ defined by sending $\phi ({\bf x} )$ to the corresponding  points in ${\bf x}$, $\psi\circ \phi$   is in  $G^0_{\Sigma, [{\bf  x}]}$. This proves that  $G_T$ and $G^0_{\Sigma, [{\bf  x}]}$ generate 	 $G_{\Sigma}$.

  To see that $G_T$ is normal,  note that for any $y\in \Sigma_v$, $\phi\in  G_{\Sigma}$  and $\psi\in G_T$, both $\phi\circ \psi (y)$ and $\psi\circ \phi (y)$ are  lying in the component
  $\phi(\Sigma_v)$. Hence the two components
  $\psi\circ\phi(\Sigma_v)$ and $\phi\circ\psi(\Sigma_v)$ are the same for any $v\in T$.
  This proves that
  $\phi\circ \psi\circ \phi^{-1}\in G_T$ and $G_T$ is normal.

  \QED

  \subsection{Moduli space of the  genus zero stable curves modeled on ${\hat T}$}
The moduli space of the  genus zero stable curves modeled on ${\hat T}$, denoted by ${\cal M}_{\hat T }$, is defined to be the equivalence classes of such curves.

Thus ${\cal M}_{\hat T }$ is the quotient space of ${\tilde {\cal M} }_{\hat T }$,  by the obvious diagonal action of $H_T=\prod_{v\in T} {\bf PSL}(2,{\bf C })_v$.
Here  ${\tilde {\cal M} }_{\hat T }$, as an open set of the product of  ${\bf CP}^1$ with $|{\bf p}|$ factors,  is the set of special points ${\bf p}=\{p_{vi}, v\in T\}$
with $p_{vi}$ being special points on $ {\bf CP}^1_v(=$ a copy of $ {\bf CP}^1$),
and  $|{\bf p}|$ is the number of spacial points.
It is easy to see that the action of $H_T$ is free and holomorphic.  One can show that the action is proper so that  ${\cal M}_{\hat T }$
is a complex  manifold. Another way to show that  ${\cal M}_{\hat T }$ is a complex manifold is to use the slice of the $H_T$-action
given by a complex submanifold of ${\tilde {\cal M}}_{\hat T }$ described in next lemma.

Next we define a global coordinate on ${\tilde {\cal M}}_{\hat T }$.
 Note that
  given a genus zero  minimal stable curve $(\Sigma,  {\bf d},{\bf x})$
  modeled on an $n$-labeled tree ${\hat T},$ there is a ``canonical'' ordering for its double points by fixing a root $v\in T$ as follows.  Using the first marked point on each top bubble component and their orders in ${\bf x}$,  all the top bubbles are ordered. Form paths of components, each  consisting of components of $\Sigma$ starting from the root component $\Sigma_v$ and connecting to one of the top bubbles. Then above ordering for the top bubbles induces an ordering for all the paths, which in turn gives an ordering of all components as well as an ordering of all double points. Here if a component
  is lying on several paths, it  belongs to the ''smallest'' path.
  This total ordering induces
  biholomorphic identifications $\phi_v:\Sigma_v\rightarrow {\bf CP}^1,$
  $ v\in T$ by sending the first three  points of ${\bf p}_v=({\bf b}_v, {\bf x}_v)$ on $\Sigma_v$ to $0, 1$ and $\infty$.

   The "global coordinate"  of  ${\cal M}_{\hat T }$ is defined by multi-cross ratios  using  the order of $p_{vi}$   above as follows.

 As mentioned above, each point in ${\tilde {\cal M}}_{\hat T }$ can be consider as a tuple of the special points on the components ${\bf CP}^1_v, v\in T.$  For each $v\in T,$ the special points $p_{vi}$ are ordered  described  above, given by the index $i$ with  $i=1, 2, \cdots, I_v.$ Then the "coordinate" for $p_{vi}$ with $3<i\leq I_v,$ is given by the  cross-ratio $w_{vi}=:w_{v123i}=(p_{v1}:p_{v2}:p_{v3}:p_{vi}).$ The tuple of all such coordinates together,  denoted by $w_{\bf p},$  with the order given by the above total ordering  is the coordinate
  for the tuple ${\bf p}=:\{p_{vi}, v\in T\}$ of special points.  The proof of the next lemma  is clear.

  \begin{lemma}
  	The map ${\bf p}\rightarrow w_{\bf p}$is $H_T$-invariant. It gives rise a
  	global  "coordinate  chart"  for  ${\cal M}_{\hat T}$. The complex submanifold of ${\tilde {\cal M}}_{\hat T}$, $S{\tilde {\cal M}}_{\hat T }=\{p_{vi}, v\in T\, |p_{v1}=0, p_{v2}=1, p_{v3}=\infty\}$ is a slice of the $H_T$-action. 
  \end{lemma}

 The universal  /tautological curve as a set is defined in the obvious way.
 One can show that  ${\cal U}_{\hat T }$  is an analytic space with only normal crossing singularities.

  \subsection{The space of  nodal maps modeled on $T$}
  Fix a Riemannian manifold $(M, g_M) $.
  Let $\Sigma$ be a nodal surface modeled on $T$.
   A  continuous map $f:\Sigma\rightarrow M$ map  is said to be a nodal $L_k^p$-map modeled on $T$
  if each its component $f_v:\Sigma_v\rightarrow M, v\in T$ is of class
  $L_k^p$. Here the $L_k^p$-norm is measured with respect the metric $g_M$ on $M$  and  the pull-back of the Fubini-Study  metric on $\Sigma_v$  by a  identification $\phi_v:\Sigma_v\simeq {\bf CP}^1.$  Note that $\phi_v$ depends on the choice of a minimal stabilization of $\Sigma $ and the total ordering of the special points,  but the induced $L_k^p$-norms are equivalent with respect to different choices.

  Two such nodal maps $f_1:\Sigma_1\rightarrow M$ and $f_2:\Sigma_2\rightarrow M$ are said to be equivalent if  is an equivalence map $\phi: \Sigma_1\rightarrow  \Sigma_2$ such that $f_1=f_2\circ \phi.$

  \begin{definition}
  	Fix a $n$-labeled tree ${\hat T}$ that is a minimal stabilization of the given tree $T$. Let ${\cal U}_{\hat T}\rightarrow {\cal M}_{\hat T}$ be the universal family of genus zero stable curves with $n$-marked points. The space of genus zero stable $L_k^p$-maps modeled on ${\hat T}$ is defined to be ${\tilde {\cal B}}^{\hat T}=$
  	$$  \{({\bf \Sigma}, f)|\, f=\{f_v:\Sigma_v\rightarrow M, v\in T\},\, f_v(d_{vu})=f_u(d_{uv })\, \, if\, \, uv\in E_T, \, \|f_v\|_{k, p}<\infty\}.$$
  	Here the domain $\Sigma$ of $f$ with components $\Sigma_v, v\in T$ is the underlying nodal curve of the stable curve ${\bf \Sigma},$ which  is  a fiber of the universal family ${\cal U}_{\hat T}\rightarrow {\cal M}_{\hat T}$.
  \end{definition}

   The space ${\tilde {\cal  B}}^{\hat T}$ can be thought as a   precise version of the intuitive notion of the space of  parametrized
    nodal maps.

  \begin{pro}
  ${\tilde {\cal B}}^{\hat T}$ is a Banach manifold of class at least $C^1$.  
  \end{pro}	

  \proof

  Recall that the universal curve can be realized as ${\cal U}_{\hat T}\rightarrow S{\cal M }_{\hat T}.$ The slice $S{\cal M }_{\hat T}$ is the product $\prod_{v\in T}{\bf S}_v,$ where ${\bf S}_v=\{p_{vi}\}$ is the set of special points on ${\bf CP}^1_v$ with  the first three of them are $0, 1, \infty.$ Note that for an unstable $v$, ${\bf S}_v$  has only one element  $(0, 1, \infty).$
  	
  	Let $${\tilde {\cal B}}_v=\{(p_v, f_v))|\, p_v\in {\bf S}_v, \, f_v:{\bf CP}^1_v\rightarrow M , \, \|f_v\|_{k, p}< \infty\}.$$ Then ${\tilde {\cal B}}_v$ is naturally identified with ${\bf S}_v\times {\tilde {\cal B}}({\bf CP}^1_v)$. Here ${\tilde {\cal B}}({\bf CP}^1_v)$ is the  Banach manifold of $L_k^p$ maps with the  fixed domain ${\bf CP}^1_v$. This implies that  ${\tilde {\cal B}}_v$ and hence
  	$\prod_{v\in T}{\tilde {\cal B}}_v$ are Banach manifolds.
  	
  	We need the following lemma.
  	\begin{lemma}
  The evaluation map $E_v:{\tilde {\cal B}}_v\rightarrow M^{|{\bf d_v}|}$given by $(p_v, f_v)\rightarrow f_v(d_v)$ is of class at least $C^1$. Here $d_v$ are the double points among the special points $p_v$ and $f_v(d_v)$  is the tuple of the evaluations of $f_v$ at the  doubles  of  $d_v$  considered   as an element in  $M^{|{\bf d_v}|}$. Moreover, the evaluation map is a surjective submersion.
  	\end{lemma}
 
\proof 

 A proof that the total evaluation  map is of class $C^{m_0}$ can be found in \cite {8}.  For the weaker result here, we  outline the proof below.
 It is easy to see that for a fixed double point $d_v$, $E_v$ is of class $C^{\infty}$. Hence for the domain  ${\tilde {\cal B}}_v\simeq {\bf S}_v\times {\tilde {\cal B}}({\bf CP}^1_v)$,  the partial derivative of $E_v$  along ${\tilde {\cal B}}({\bf CP}^1_v)$-direction is at least of class $C^\infty$. By our assumption and Sobolev embedding theorem, the partial derivative of $E_v$  along ${\bf S}_v$ is of class $[m_0]\geq 1.$
  
  \QED

  Let $E=\prod_{v\in T}E_v:{\tilde {\cal B}}_v\rightarrow M^{|{\bf d}|}$
  be the total evaluation map at double points. Then  it is still a submersion so that it is  transversal to the diagonal $\Delta_T\subset  M^{|{\bf d}|}.$ Here $\Delta_T=\{ m_{uv}=m_{vu}\in M, when \, \,uv \in E_T\, \,for \, \, u, v \in T  \}$ is a submanifold of $  M^{|{\bf d}|}=\{m_{uv}\in M, uv\in E_T\}.$  This implies that ${\tilde {\bf B}}^{\hat T}=E^{-1}(\Delta_T)$ is a Banach manifold of class at least  $C^1.$

  \QED

In Gromov-Witten theory we are mainly interested
in the space of equivalent classes of the nodal $L_k^p$-maps, denoted
by ${\cal B}^T.$  It can be obtained by quotient out the action of the reparametrization groups on
${\tilde {\bf B}}^{\hat T}$.   These actions are induced from the corresponding  actions on the domains: $G_{\Sigma}\times \Sigma\rightarrow \Sigma$ for ${\bf \Sigma}\in  {\cal M}_{\hat T}$.
There is a coherency of these actions as ${\bf \Sigma}$ moving in $  {\cal M}_{\hat T}$ as follows.

Let $N=:N_{\epsilon}({\bf \Sigma}_0)$ be a small neighborhood  of a stable curve
${\bf \Sigma}_0=(\Sigma_0, {\bf d}_0, {\bf x}_0)$ with $n$ marked points in ${\cal M }_{\hat T}$, and  ${\cal U}_{\hat T, N}\rightarrow N$ be the local  universal family. Identify $N$ with a poly-disk with coordinate $t$ with $t=0$ for ${\bf \Sigma}_0$.  Let  $(\Sigma_0, {\bf x}_0)$ be the central fiber $({\cal U}_{\hat T, N})_0$ of the local  universal family.

\begin{pro}
The actions of $G_{\Sigma_0}$ and $G_{\Sigma_0, [x_0]}$ on the central fiber extends as continuous actions of the local universal family ${\cal U}_{\hat T, N}\rightarrow N_{\epsilon}({\bf \Sigma}_0)$
when $\epsilon$ small enough. For any ${\bf \Sigma}_1\in N_{\epsilon}({\bf \Sigma}_0)$, the actions $G_{\Sigma_1}$ and  $G_{\Sigma_1, [{\bf x}_1]}$ are induced from the extended actions of $G_{\Sigma_0}$ and  $G_{\Sigma_0, [{\bf x}_0]}$ respectively. Moreover  any two fibers
${\bf \Sigma}_1$ and ${\bf \Sigma}_2$  are equivalent, considered as nodal curves and unordered marked curves respectively if and only if they are in the same $G_{\Sigma_0}$-orbit and  $G_{\Sigma_0, [{\bf x}_0]}$-orbit accordingly.

\end{pro}
 
 The part of the  above statement on  $G_{\Sigma_0}$ can be derived form the corresponding one on $G_{\Sigma_0, [x_0]}$ by noting that $G_{\Sigma_0}$ is the semi-direct product of $G_T$ and $G_{\Sigma_0, [x_0]}$.
On the other hand, the statement on $G_{\Sigma_0, [x_0]}$  is  the so called  universal property of the local  universal family, which  is  part of  the standard theory of the Deligne-Mumford type of moduli spaces. We  refer the readers to \cite{3} for a discussion on this  and the reference therein. 
Note that the action of $G_T$ is globally defined on the universal family  ${\cal U}_{\hat T}\rightarrow {\cal M}_{\hat T}$ and acts on ${\cal M}_{\hat T}$ as the
identity map. 

Thus these local actions
 together form a functorial system. So do the induced actions
 on ${\tilde {\bf B}}^{\hat T}$. They are so called groupoids (or Artin stacks in the setting of algebraic geometry). The quotient space of   ${\tilde {\bf B}}^{\hat T}$ by these actions, denoted
 by ${\cal B}^T$, is the space of  unpapametrized nodal $L_k^p$-maps.

\begin{pro}
	The action map ${\Phi}: G_T\times {\tilde {\cal  B}}^{\hat T}\rightarrow 	 {\tilde {\cal  B}}^{\hat T}$ is continuous. 
 The same is true  for the actions of  $G_{\Sigma}$.
\end{pro} 	

This  proposition follows from the well-known fact that the action map on the mapping space ${\cal M}_{k,p}({\bf P}^1, M)$ by reparametrizations is continuous. A proof can be found  for instance in \cite {5}.

 \section {The properness of $G_T$-action and \\ $G_T$-Hausdorffness of ${\cal B}^{\hat T}$}

In  this  section, for each $f\in   {\tilde {\cal  B}}^{ws}_{\hat T}$ we will consider the induced action of  $G_f$ on  the $G_f$-orbit of a prescribed neighborhood $U_f\subset {\tilde {\cal  B}}^{ws}_{\hat T}$ of $f$. Here and throughout this section we use $G_f$ to denote $G_{\Sigma_f}$  with   ${\Sigma_f}$  being the domain of $f$, and $G_{\Sigma_f}$ was defined in last section. As noted in  the previous  section, the collection of such groups together with the induced actions above form a functorial system. 

First let us prove that  for  locally compact $X$ and $G$, the definition of properness introduced in Sec. 1 is   equivalent to the usual one.

\begin{lemma}
	For  locally compact $X$ and $G$, a group action $\Phi:G\times X\rightarrow X$ is proper in the sense of this paper is 	equivalent to the usual one.
	
\end{lemma}

\proof

Assume the usual definition of the properness for $\Phi$. Let $K\subset  X\times X$ be a compact subset of
$ X\times X$. Let $U$ be an open neighborhood of $K$ obtained from a finite open covering
 of $K$ with each open set  having compact closure so that the closure ${\overline U}$ of $U$ is compact. Then   $(\Phi\times{\it id}_{X})^{-1}({\overline U})$ is compact by the assumption.
Hence  $\pi_G ((\Phi\times{\it id}_{X})^{-1}({\overline U}))$ is  compact in $G$. Since $\pi_G ((\Phi\times{\it id}_{X})^{-1}({ U}))\subset \pi_G ((\Phi\times{\it id}_{X})^{-1}({\overline U}))$, it is contained  in a compact set in $G$.

Conversely assume the definition of the properness for $\Phi$ used in this paper.
Let $K\subset  X\times X$ be a compact subset of
$ X\times X$ and  $U$ be the corresponding  open neighborhood of $K$ such that
 $\pi_G ((\Phi\times{\it id}_{X})^{-1}({ U}))$ is contained  in a compact set $K_G$ in $G$. Since $X\times X$ is locally compact, there is a neighborhood $U_1$ of $K$ with  compact  closure ${\overline U_1}$ containing inside $U$.
 Hence  $\pi_G ((\Phi\times{\it id}_{X})^{-1}({\overline  U_1}))$ is contained  in the  compact set $K_G$ as well.

Now $(\Phi\times{\it id}_{X})^{-1}(K)\subset (\Phi\times{\it id}_{X})^{-1}({\overline  U_1})\subset \pi_G ((\Phi\times{\it id}_{X})^{-1}({\overline  U_1}))\times K \subset K_G \times K.$
  Thus  $(\Phi\times{\it id}_{X})^{-1}(K)$ is a closed subset of the compact set
$K_G \times K$ in $G\times X$, hence compact.

\QED

\begin{lemma}

	A group action $\Phi:G\times X\rightarrow X$ is proper	if and only if for any compact subset $K_1\times K_2\subset X\times X$ there is a neighborhood $U_1\times U_2$ of $K_1\times K_2$ in  $ X\times X$ and a  compact subset $C$ of $G$ accordingly such that for any $h\in G\setminus C$, $ h(U_2)\cap  U_1=\varphi$.
\end{lemma}

\proof

  Note that $\pi_G ((\Phi\times{\it id}_{X})^{-1}(U_1\times U_2))=\{g\in G|\, g(U_2)\cap  U_1\not =\varphi\}.$
  Hence if $\pi_G ((\Phi\times{\it id}_{X})^{-1}(U_1\times U_2))$ is pre-compact, let $C$ be the closure of  $\pi_G ((\Phi\times{\it id}_{X})^{-1}(U_1\times U_2))$. Then for any $h\in G\setminus C$, $ h(U_2)\cap  U_1=\varphi$.

 Conversely if there is a compact subset $C$ of $G$ such that for any $h\in G\setminus C$, $ h(U_2)\cap  U_1=\varphi$, then $h\in G\setminus \pi_G ((\Phi\times{\it id}_{X})^{-1}(U_1\times U_2 ))$. Hence  $\pi_G ((\Phi\times{\it id}_{X})^{-1}(U_1\times U_2 ))$ is contained in the compact set $C$.

 \QED

\begin{lemma}
	A group action $\Phi:G\times X\rightarrow X$ is proper	if and only if for
	any point $p\in X\times X$ there is a neighborhood $U$ of $p$ in  $ X\times X$
	such that
	the closure of  $\pi_G ((\Phi\times{\it id}_{X})^{-1}(U))$ is compact in $G$.

\end{lemma}
	
\proof

 One direction is clear. For the other direction,
given a compact set $K\subset  X\times X$,
for each point
$p\in K $, let $U_{p}$ be the  corresponding
neighborhood given by the  assumption of the lemma so that $\pi_G ((\Phi\times{\it id}_{X})^{-1}(U_p ))$ is contained in a compact set. Then there are finitely many such points in $K$, denoted by $p^i, $ $i=1, \cdots, I$  and the corresponding
neighborhoods, denoted by   $U(p^i)$ such that
the open set  $U_K$ defined by $U_K=\cup_{i=1}^{ I}U(p^i)$ contains $K$. Then $\pi_G ((\Phi\times{\it id}_{X})^{-1}(U_K ))=\cup_{i}^{ I}\pi_G ((\Phi\times{\it id}_{X})^{-1}U(p^i)).$
 Since each $\pi_G ((\Phi\times{\it id}_{X})^{-1}U(p^i))$ is contained in a compact set
in $G$ by our assumption, so is $\pi_G ((\Phi\times{\it id}_{X})^{-1}(U_K ))$.

\QED

By the above lemma,  the properness of the action of $G_f$ on  $   {\tilde {\cal  B}}^{ws}_{\hat T}$ can be derived from the following theorem.

 \begin{theorem}
 	The partially defined action of $G_f=:G_{\Sigma_f}$ on  $   {\tilde {\cal  B}}^{ws}_{\hat T}$ has the following property:
 	 for any $f_1$ and $f_2$  $   {\tilde {\cal  B}}^{ws}_{\hat T}$ there exist the   open neighborhoods $U_1$ and $U_2$ containing $f_1$ and $f_2$ with
 	$G_{f_1}$ and $G_{f_2}$ acting upon  respectively  and compact subsets $K_1$ and $K_2$ in 	$G_{f_1}$ and $G_{f_2}$ accordingly   such that  for any $h_1$ in $U_1$ ($h_2$ in $U_2$) and $g_1$ in $G_{f_1}\setminus K_1$  ($g_2$ in $G_{f_2}\setminus K_2$ ), $g_1\cdot h_1$ is not in $U_2$ ($g_2\cdot h_2$ is not in $U_1$).

 \end{theorem}

\proof

The  theorem  can be derived from the corresponding main theorem on properness of $G_T$-action in section 1.

  Indeed,  if the theorem is not true, we may assume that there  are  sequences of $h_i\in
  U_1$ with $\lim _{i\rightarrow \infty}h_i=f_1$ and $g_i\in G_{f_1}$ with $g_i$ not staying in any compact set as $i$ goes to infinity such that $g_i\cdot h_i=h_i\circ g_i$ is lying in $U_2$ with $\lim _{i\rightarrow \infty}h_i\circ g_i=f_2$. Now
  $g_i=p_i\circ t_i$ with $t_i$ in $G_T$ and $p_i$ in  the finite group
  $G_{\Sigma_{f_1, [{\bf x}]}}.$ Hence after taking a subsequence, we may assume that $p_i$  is a fixed element $p$ in $G_{\Sigma_{f_1, [{\bf x}]}}$.
   Rename $f_1\circ p$ by $f_1$,  $p(U_1)$ by $U_1$ and $h_i\circ p$ by $h_i.$
  
   Then we get  the sequences of $h_i\in
   U_1$ with $\lim _{i\rightarrow \infty}h_i=f_1$ and $g_i\in G_{T}$ with $g_i$ not staying in any compact set as $i$ goes to infinity such that $g_i\cdot h_i=h_i\circ g_i$ is lying in $U_2$ with $\lim _{i\rightarrow \infty}h_i\circ g_i=f_2$.
  This contradicts to the main theorem
  on the properness of $G_T$-action.

  Clearly to prove the main theorem for $G_T$, we only need to look its action on its unstable components. Since $G_T$-action on the domains maps each component to itself and  does not move the  double points (if there is any), we only need to prove the corresponding statement for a fixed domain $\Sigma\simeq {\bf CP}^1$
  with $G_T=G_i, i=0, 1, 2$ of the subgroups of ${\bf PSL}(2, {\bf C})$ preserving $i$ marked points. The  proofs of the three cases are similar. We  give the full detail of the proof for  the hardest case that  $G_T=G_0= {\bf PSL}(2, {\bf C})$ in the next proposition. Clearly the case $G_T=G_2= {\bf C}^*$ is a special case of this.   As for $G_T=G_1$, a detailed proof  was given in \cite{6}. For completeness, we recall the argument there in
the  proposition after  next one.

\QED
  \begin{pro}
  	Let $   {\tilde {\cal  B}}^{*}$ be the Banach manifold of non-trivial $L_k^p$-maps from $\Sigma={\bf CP}^1$ to $M$. Then the  action of  $G={\bf PSL}(2, {\bf C})$ on $   {\tilde {\cal  B}}^{*}$ is proper.
  \end{pro}

  \proof

 We start with some elementary linear algebra.

 For any  $g \in SL(2, {\bf C}),$
 let    $g=h\cdot  u$ be the "polar decomposition"   with $u\in U(2)$ and $h$ being self-adjoint and positive.  In fact here   $h=(g\cdot g^*)^{\frac {1}{2}}$ is in $ SL(2, {\bf C})$  so that  $u= (g\cdot g^*)^{-\frac {1}{2}}\cdot g\in SU(2)$.

  Then  $(g\cdot g^*)^{\frac {1}{2}}=w^*\cdot diag ( r_1, r_2)\cdot w$  with $w\in SU(2)$ and  $r_1\leq r_2$.   Rename $w^*$ as $u$ and $wu$ as $v$. Denote $diag ( r_1, r_2)$ by $\Delta({\bf r})$ for short. Then we have a decomposition  $g=u\cdot \Delta({\bf r})\cdot v$ in $SL(2,{\bf C})$ with $u$ and $v$ in
 $SU(2)$. Since $g\cdot g^*=u\cdot \Delta({\bf r})^2\cdot u^*$, the decomposition is essentially unique if $g\cdot g^*\not= I_2$ through we do not need this.

  Assume that the proposition  is not true. There are two possibilities : (i) for any small neighbourhoods
  $U_{\epsilon_i}(f_i), i=1, 2 $ and a compact subset $K_n$ in $G$,
  there is a $g_n\in G$ and $h_n\in U_{\epsilon_1}(f_1)$ such that (a) $g_n$ is not in $K_n$; (b) $h_n\circ
  g_n$ is in  $U_{\epsilon_2}(f_2);$ (ii)  the same statement as (i) but switching the roles of $f_1$ and $f_2$.
  The two case are symmetrical. Hence we may  assume that we are in case (i).

   \vspace{2mm}
   \noindent

 For this case,  we  fix any  $U_{\epsilon_1}(f_1)$,    and construct $U_{\epsilon_2}(f_2)$ and a  nested infinite  sequence of compact sets $K_1\subset K_2\subset \cdots \subset K_l\cdots$ in $G$ below. Assuming this is done,  then the assumptions above imply there exist the corresponding  sequences
  $\{g_n\}_{n=1}^{\infty}$ and $\{h_n\}_{n=1}^{\infty}$ with the properties (a) and (b)
 above. We will show then  (a) and (b) lead to a contradiction.

 To construct $K_n$, let  ${\tilde K}_n\subset  U (2)\times {\cal D}_{2} \times U (2)$   be the compact subset of tuples  $(u, D({\bf r}), v)$  with $\frac {1}{n} \leq r_i\leq 1, i=1, 2$, where $\Delta({\bf r})=diag (r_1,
 r_{2})$ and ${\cal D}_{2}$  is the collection of all $2\times 2$
 non-singular  diagonal matrices  with positive  entries, and $\pi: U (2)\times {\cal D}_{2} \times U (2)\rightarrow {\bf PSL}(2, {\bf C})$ given by mapping $(u, \Delta({\bf r}), v)$   to the image of $u\cdot \Delta({\bf r})\cdot v$ in  ${\bf PSL}(2, {\bf C})$.  Denote  the corresponding compact set in  ${\bf PSL}(2, {\bf C})$ by ${ K}_n$. Note that the image of  $U(2)$  under the natural projection
  to  $GL(2, {\bf C})/{\bf C}^*\simeq {\bf PGL}(2, {\bf C})={\bf PSL}(2, {\bf C})$
 is contained in all ${ K}_n$.

 Next we choose $\epsilon_2$ for $U_{\epsilon_2}(f_2)$  as follows.    Since $f_2$ is nontrivial, its  energy $E(f_2)=\delta_2>0.$ Then there exits a point $x_0\in {\bf CP}^1 $ such that $e(df_2)(x_0)=:|df(x_0)|^2>0.$  Hence  there are   positive constants $\gamma$ and $\rho$ small enough such that for any
 $x$ in the disc  $B(x_0; \rho)$  of radius $\rho$ centered at $x_0$,  $e(df_2)(x)>\gamma$. Then there is an ${\tilde \epsilon}_2>0$ such that for any $h$ with $\|h-f_2\|_{C^1}<{\tilde  \epsilon}_2$, the same is true that
 $e(dh)(x)>\gamma$ for  $x\in B(x_0; \rho)$.

 Now $\|h-f_2\|_{C^1}\leq C_2 \cdot \|h-f_2\|_{k, p}$ by our assumption.
 Hence  we choose $\epsilon_2$ by the requirement that $\epsilon_2<{\tilde \epsilon}_2/C_2$ so that  for any $h\in U_{\epsilon_2}(f_2)$,  $\|h-f_2\|_{C^1}
 <{\tilde \epsilon}_2$.

With this choice of $\epsilon_2$, for any $h\in U_{\epsilon_2}(f_2)$, and any point $x_1\in B(x_0; \rho)$ and  $\rho_1<<  \rho$ such that  $B(x_1; \rho_1)\subset B(x_0; \rho),$  the  energy $E(h|_{ B(x_1; \rho_1)})>\pi\gamma\cdot \rho^2_1$.

  By our definition of $K_n$, we may assume that for the lifting $g_n$ in ${ SL}(2, {\bf C}), $ denoted by the same notation, has a  decomposition in ${\bf SL}(2, {\bf C}), $
  $g_n=u_n\cdot \Delta({\bf r}_n)\cdot v_n$ with $u_n$ and $v_n$ in $SU (2)$ and
  $ \Delta({\bf r}_n)=diag (r_{n,1}, r_{n,2})$ with  $0<r_{n,1}\leq  r_{n,2}.$ We may assume that $ r_{n,2}=1$ by considering $ \Delta({\bf r}_n)$ as an element
  in ${\bf PGL}(2, {\bf C}).$
  Denote  $r_{n,1}$  by $a_n$  and
  $ \Delta({\bf r}_n)$ by  $\Delta(a_n)$.

  In these notations, the condition (a) above implies that for $g_n=u_n\cdot \Delta({a}_n)\cdot v_n$,     $lim_{n\mapsto \infty} a_n=0.$
After  taking subsequence, we may assume that $lim_{n\mapsto \infty} u_n=u$ and
  $lim_{n\mapsto \infty} v_n=v$  in $SU(2)$.   Note that when considered as  automorphisms on ${\bf CP}^1$, the convergence  here are with respect to $C^{\infty}$-topology on the  corresponding mapping space.

 Let  ${\bf CP}^1={\bf C}\cup {\infty}  $. Now we need deal with two case:
 (I)  $\infty\in v(B(x_0; \rho) $ and (II) $\infty \not \in v(B(x_0; \rho) $.

  Since $v$ is an isometry, it is easy to see that in both cases, there exits an $x_1\in B(x_0; \rho) $ and a positive number $\rho_1<< \rho$ such that $\infty$  is not lying in the closure of $v(B(x_1; \rho_1)).$
We may assume that $dist (v(B(x_1; \rho_1)),\infty)>\delta>0$. Then 
for $i>i_0$ large enough, $dist (v_i(B(x_1; \rho_1)),\infty)>\delta>0$ as well.
Hence
there is a large $R$ such that  for $i>i_0$ large enough, $v_i(B(x_1; \rho_1))$ is lying inside 
$D_R\subset {\bf C}^n$   of the  open disks centered at origin  in ${\bf C}$  with  radius $R$. Fix such a $R$.

 Note that  in term of the   coordinate  of  ${\bf C}\subset {\bf CP}^1$, the action of $\Delta({a_i})$  is given by $\Delta{a_i}(z)= a_i\cdot z.$
 Hence for any fixed $R>0 $ and any given $\epsilon >0, $    our assumption that $a_i\mapsto 0$ implies that there is a fixed  $i_0(\epsilon)>>0$ such that when $i>i_0(\epsilon),$
 $  \Delta({a_i})(v_i(B(x_1; \rho_1)))=a_i\cdot (v_i(B(x_1; \rho_1))\subset a_i\cdot B_R\subset B_\epsilon\subset  {\bf C}$.
Hence the areas $A (\Delta({a_i})\circ v_i( B(x_1; \rho_1)) )\leq C_3\epsilon^2.$   Here the areas are computed with respect to the Fubini-Study metric which is uniformly equivalent
 to the flat metric on $B_R$ for fixed $R$.

 Applying this to $g_i=u_i\circ \Delta({a_i})\circ v_i$, since $u_i$  preserves the  Fubini-Study metric, we conclude that
 for $i$ large enough,   $A (g_i(B(x_1; \rho_1)))<C_3\epsilon^2.$

 Now

 $$E  ((h_i\circ g_i)|_{ B(x_1; \rho_1)})=E  (h_i|_{ g_i( B(x_1; \rho_1))})$$  $$ \leq ||h_i||^2_{C^1} A ( g_i( B(x_1; \rho_1)))\leq C_3\cdot ||h_i||^2_{C^1}\epsilon^2$$
 $$  \leq C_3 ||h_i||^2_{k, p} \epsilon^2  \leq  C_3(\|f_1\|_{k, p}+ ||h_i-f_1||_{k, p})^2\epsilon^2 $$ $$ \leq C_3(\|f_1\|_{k, p}+\epsilon_1)^2 \epsilon^2 . $$
 By letting $\epsilon \mapsto 0,$ we conclude that  $\lim_{i\mapsto\infty}E  ((h_i\circ g_i)|_{ B(x_1; \rho_1)})
 =0.$ 

 Now  $h_i\circ g_i\in U_{\epsilon_2}(f_2)$, and  recall that  for any $h\in U_{\epsilon_2}(f_2)$,  the  energy $E(h|_{ B(x_1; \rho_1)})>\pi\gamma\cdot \rho^2_1>0$.
  This is a contradiction.

  \QED

 \begin{pro}
 The  action of  $G_1 $ on $   {\tilde {\cal  B}}^{*}$ is proper.
 \end{pro}

 \proof

Recall that in this case, $G=G_1$ is  the automorphisms of  $ ({\bf CP}^1, d)=({\bf C}\cup \{{\infty}\}, \infty)$ with a fixed point $d=\infty$.
In the coordinate $z\in {\bf C}\subset {\bf CP}^1, $
 each element $g\in G$ has the form $g(z)=a(z-c)$ with $a\not =0.$
Since $f_1$ is not a constant map,
we may assume that, $f_1(0)\not= f_2(\infty)$ by choosing a different origin of
${\bf C}$.
Assume that  $d(f_1(0),  f_2(\infty))=\delta_0>0$.

The first requirement on  $\epsilon_i$.
 is  that
$\epsilon_i<< \delta_0, i=1,2.$

Since  $f_2$ is not a constant map, there is $x_0\in {\bf C}$ and $\rho>0$ and $\gamma>0$ such that $e(df_2)(x)>\gamma$ for any $x\in D(x_0, \rho)\subset {\bf C}.$
The second requirement for $\epsilon_2$ is that  for any $h\in U_{\epsilon_2}(f_2)$ the same is true.

Now choose  a sequence of compact subsets ${\tilde K}_n\subset C^*\times {\bf C}$  defined by  ${\tilde K}_n=\{(a, c)\in C^*\times {\bf C}\,|\,\frac {1}{n} \leq |a|\leq n, \, |c|\leq n\}$. Let $K_n $ be the corresponding compact subsets in $G$ given by the map $C^*\times {\bf C}\rightarrow G$ sending $(a, c)$ to $g(z)=a(z-c)$.

Assume that the Theorem  is not true for $G_1$ so that  for all $n$ there  are
$g_n(z)=a_n(z-c_n)$ not in $K_n$ and $h_n$ in $U_{\epsilon_1}(f_1)$ such that  $h_n\circ
g_n$ is in  $U_{\epsilon_2}(f_2).$

\vspace{2mm}
\noindent
\vspace{2mm}
\noindent
$\bullet$   Claim:
$|c_n|$ is bounded.

\proof

Assume that the claim is not true. Then there is a subsequence of $\{g_n\}_{n=1}^{\infty}$, denoted by the same notation, such that  $c_n\rightarrow \infty$  as $n\rightarrow \infty.$

Now  $h_n\circ g_n(c_n)=h_n(a_n(c_n-c_n))=h_n(0).$ Since $h_n\in U_{\epsilon_1}(f_1),$
we have $$d(h_n\circ g_n(c_n), f_1(0))=d(h_n(0), f_1(0))$$ $$\leq \|h_n-f_1\|_{C^0}
\leq C_1\|h_n-f_1\|_{k,  p}\leq C_1\epsilon_1.$$

Since $h_n\circ g_n$ is in $ U_{\epsilon_2}(f_2),$  we have that $$d(h_n\circ g_n(c_n), f_2(c_n))\leq \|h_n\circ g_n-f_2\|_{C^0}
\leq C_2 \|h_n\circ g_n-f_2\|_{k,  p}\leq C_2\epsilon_2.$$

The continuity of $f_2$ at $\infty$ implies that for any given  $\epsilon$,
$d(f_2(c_n), f_2(\infty))\leq \epsilon$ when
$n$ is large enough. Therefore,  $d(h_n\circ g_n(c_n), f_2(\infty))\leq  C_2\epsilon_2+ \epsilon .$       We conclude that $\delta_0=d(  f_1(0),   f_2(\infty))
\leq C_1\epsilon_1+ C_2\epsilon_2+ \epsilon$. This contradicts to our assumption that $\delta_0>>\epsilon_i, i=1,2.$

\QED

Therefore, we may assume that $\lim_{n\mapsto\infty}c_n=c \in {\bf C}.$
By shifting the origin again, we may assume that $c=0.$ With  respect to this new
 coordinate $z\in {\bf C}\subset {\bf CP}^1, $
  we still have  $g_n(z)=a_n(z-c_n)$ but $\lim_{n\mapsto \infty}c_n=0.$
Then our assumption  implies that
the sequence $\{a_n\}_{n=1}^{\infty}$ in $C^*$   either going  to zero or to infinity. 

In this new coordinate, $D(x_0, \rho)$ is still lying in ${\bf C}$.
By replacing $D(x_0, \rho)$ by a sub-disk $D(x_1, \rho_1)\subset D(x_0, \rho)$, we may assume that $0\not \in D(x_0, \rho).$
Then
 there  exist  $0<R_1<R$,  such that
$D(x_0, \rho)$ is contained inside the annuals $B_R\setminus B_{R_1}$.
For $n$ sufficiently large, since $c_n\mapsto 0$ as $n\mapsto \infty,$ we may assume that the sets $D(x_0, \rho)\pm c_n$  are still contained inside the annuals $B_R\setminus B_{R_1}$.

 We now  assume that $\lim_{n\mapsto \infty} a_n=0.$  The other case  that $\lim_{n\mapsto \infty} a_n=\infty$
 can be treated symmetrically by  the  construction above of $D(x_0, \rho)$ and $B_R\setminus B_{R_1}$.

 Then for $n$ sufficiently large,
$$E(h_n\circ g_n|_{D(x_0, \rho)})= E(h_n |_{g_n(D(x_0, \rho))})\leq  E(h_n |_{a_n(B_{R})})\leq C\cdot a_n^2A(B_R)(\|f_1\|_{k, p}+\epsilon_1)^2.$$

 Hence $E(h_n\circ g_n|_{D(x_0, \rho)})$ tends to zero as $n$ tends to infinity.

On the other hand since  $h_n\circ g_n$ is inside $U_{\epsilon_2}(f_2)$, by the choice of $\epsilon_2$, $E(h_n\circ g_n|_{D(x_0, \rho)})>\pi\rho^2\gamma.$
This is a contradiction.
\QED

Applying  the theorem on proper action  to the case  that $f_1=f_2$,  we get the first  part of the following  corollary.

 \begin{cor}
 	The isotropy group ${\Gamma}_f$ of any weakly stable $L_k^p$-map
 	$f$ is always compact. Moreover the  $G_T$-orbit and hence $G_f$-orbit of $f$  in
 	$   {\tilde {\cal  B}}^{\hat T}$ is closed.

 \end{cor}

 \proof

 We only  need to  prove  the  last statement.
 We only need prove this   for the essential case that $f\in  {\tilde {\cal B}}^{ws}_{\hat T}$.
 Rename $f$ as  $f_1$. If the corollary is not true, there exist $g_i\in G$ and $f_2\in {\tilde {\cal B}}^{ws}_{\hat T}$ such that
 $f_2=\lim_{i\mapsto \infty}f_1\circ g_i,$ but $f_2$ is not in $G\cdot f_1$.
 Therefore for any $U_{\epsilon_2}(f_2), $  when $i$ is large enough,  $f_1\circ g_i$ is in $U_{\epsilon_2}(f_2). $
 On the other hand,  the main theorem  with the same notation implies that for all such $i$,  $g_i$ is in the compact set $K_1$. Therefore, we may assume that
 $\lim_{i\mapsto \infty} g_i=g$ in $K_1.$ Consequently, $f_2=\lim_{i\mapsto \infty}f_1\circ g_i=f_1\circ g.$ That is $f_2\in G\cdot f_1$ which is  a contradiction.

 \QED

 Next theorem is a slightly stronger version of the Theorem 1.1

 Recall that  given  $f$ in   $   {\tilde {\cal  B}}_{\hat T}$, there is sufficiently small neighborhood $U=U_f$ such that the action $G_f$ extends to
 $U$ so that for  any $h\in U$ the action of $G_f$ covers the action
 of $G_h$ on the  corresponding $U_h\subset U_f$.

  \begin{theorem}
  The space  $   {\tilde {\cal  B}}^{\hat T}$ of  nodal $L_k^p$-maps modeled on ${\hat T}$  is $G_f$-Hausdorff in the sense that for any two diffent $G_f$-orbits $G_{f_1}f_1$ and $G_{f_2}f_2$, there exit $G_f$-neighborhoods $G_{f_1}U_1$ and $G_{f_2}U_2$ such that $G_{f_1}U_1 \cap G_{f_2}U_2=\varphi$. Therefore, the  quotient space $  { {\cal  B}}^{\hat T}$  of unparametrized  nodal $L_k^p$-maps is Hausdroff.
  \end{theorem}

 \proof

We make a few  reductions.

By the definition of the  partially defined actions of $G_{f_i}, i=1, 2$ on $   {\tilde {\cal  B}}^{\hat T}$,  the actions are the ones induced from the  action on the special points of the domains, considered as elements in the slice  $S{\cal M }_{\hat T}=\prod_{v\in T}{\bf S}_v,$ where ${\bf S}_v=\{p_{vi}\} $ is the set of special points on ${\bf CP}^1_v\simeq {\bf CP}^1$ with  the first three of them are $0, 1, \infty$ (see the definition in Sec. 2 ).
Note  $G_{f_i}=G_{{\Sigma}_i}$ is the semi-direct product of $G_T$ and $G^0_{{\Sigma}_i, [{\bf x}_i]}$.
 
 Let $(\Sigma_1, {\bf d}_1, {\bf x}_1)$ and $(\Sigma_2, {\bf d}_2, {\bf x}_2)$ be  minimal stabilizations of the domains for $f_1$ and $f_2$.
  Assume that for these  domains, the $G_{f_1}$ ($=G_{{\Sigma}_1}$)-orbit of  $(\Sigma_1, {\bf d}_1, [{\bf x}]_1)$ and $G_{f_2}$ ($=G_{\Sigma_2}$)-orbit of $(\Sigma_2, {\bf d}_2, [{\bf x}]_2)$  do not intersect when considered as nodal surfaces.  
    Then assumption here  simply means  that the  two orbits for the corresponding sets of the double points on the stable components   do not intersect.
  Since only the finite subgroup $G_{{\Sigma}_i, [{\bf x}_i]}$, $i=1, 2$ act on 
  the double points on the stable components,  we conclude  that for sufficiently
   small neighborhoods $N_{\epsilon_i} (\Sigma_i, {\bf d}_i, {\bf x}_i)$  the corresponding two $G_{f_i}$-orbits  do not intersect as well.
   By definition, the theorem holds in this case simply because the domains of $f_1$ and $f_2$ are already separated by $G_{{\Sigma}_i, [{\bf x}_i]}
   	$-neighborhoods.
   
   Thus we may assume that the orbits for  domains   intersect so that we only need consider the case that they have the same domains, though in the following we still use different notations to denote the domains and related objects.

First note that if $G_{f_1}U_1 \cap G_{f_2}U_2\not =\varphi$ implies that there exist $\phi_1\in {\Sigma_{f_1, [{\bf x}_1]}}$ and $\phi_2\in {\Sigma_{f_1, [{\bf x}_2]}}$ such that $G_{T_1}\phi_1(U_1) \cap G_{T_2}\phi_1(U_2)\not =\varphi$. Since $f_1\circ \phi_1$ and $f_2\circ \phi_2$
satisfy the  corresponding condition  in the theorem,  we may replace $f_1$ and $f_2$ by $f_1\circ \phi_1$ and $f_2\circ \phi_2$ and prove the theorem for $G_T$ only.

 Next note that the condition of the theorem implies that there are  two  components, one on $f_1$, the other on $f_2$, defined on the same domain $\Sigma_v\simeq S^2$ with $v\in T$ such that the corresponding orbits do not intersect. .

 Thus we only need consider  the  case that $f_1$ and $f_2$ only have  one component   with the same domain  $\Sigma\simeq S^2$  and $G_T$-action.   The proof for this case consists of three parts.

  ${\bullet}$   Part I:  the case that   both $f_1$ and $f_2$  are  non-trivial.

  By Theorem 3.1, for any $g$ not in the compact set $ K_1$ and $h\in  U_{\epsilon_1}(f_1)$, $h\circ g$ is not in $U_{\epsilon_2}(f_2)$. By our assumption, we may assume that
  $  U_{\epsilon_1}(f_1)$ and $U_{\epsilon_2}(f_2)$ have no intersection.

  \vspace{2mm}
  \noindent
  \vspace{2mm}
  \noindent $\bullet$ Claim: when $\epsilon_i,i=1,2$  are small enough,
  $(G_T\cdot U_{\epsilon_1}(f_1))\cap U_{\epsilon_2}(f_2)$ is empty.

  \proof

  \vspace{2mm}
  \noindent
  If this is not true, there  are $h_i\in U_{\delta_i}(f_1)$ and $g_i\in K_1$
  such that $h_i\circ g_i$ is in $ U_{\delta_i}(f_2)$ with $\delta_i\mapsto 0.$
  The compactness of $K_1$ implies that after taking a subsequence, we have that
  $\lim_{i\mapsto \infty}g_i=g\in K_1.$ Since $\delta_i\mapsto 0$, we have that
  $f_1=\lim_{i\mapsto \infty}h_i$ and $f_2=\lim_{i\mapsto \infty}h_i\circ g_i=f_1\circ g.$
  Hence, $f_1$ and $f_2$ are in the same orbit which contradicts to our assumption.
  Note that in the last identity above, we have used the fact that the action map
  $\Psi: G_T\times   {\tilde {\cal  B}}^{\hat T}\rightarrow   {\tilde {\cal  B}}^{\hat T}$ is continuous.

  \QED

  Of course the same proof also implies that
  $(G_T\cdot U_{\epsilon_2}(f_2))\cap U_{\epsilon_1}(f_1)$ is also empty for sufficiently
  small $\epsilon_i, i=1, 2.$

  If $h\in (G_T\cdot  U_{\epsilon_1}(f_1))\cap ( G_T\cdot  U_{\epsilon_2}(f_2)), $  then there are
  $h_i\in U_{\epsilon_i}(f_i) $ and $g_i\in G, i=1, 2$ such that
  $h=h_1\circ g_1=h_2\circ g_2. $ Hence $h_2=h_1\circ g_1\circ g_2^{-1}$ and
  $(G_TU_{\epsilon_1}(f_1))\cap U_{\epsilon_2}(f_2)$ is not empty. This contradicts
  to the above claim.

  \QED

 ${\bullet}$   Part II:  the case that   one of  $f_1$   and  $f_2$  is
  trivial  but the other   is not.  Then the desired results follows from the following stronger statement.

  \begin{lemma}
  	Given any two $L_k^p$-maps $f_1$ and $f_2$ with energy $E(f_1)\not =E(f_2)$, there exit   $G$-neighborhoods  $  { W}(f_1)$ of $f_1$ and $ { W}(f_2)$ of $f_2$  which  do  not intersect. 
  	In particular if $f_1$ is a constant map and $f_2$ is not, then $E(f_1)=0\not =E(f_2)$ and the above conclusion  holds.
  	
  \end{lemma}

  \proof

  Note that the condition $E(f_1)\not =E(f_2)$ implies that $f_1$ and $f_2$ are not in the same $G$-orbit.
  We may assume that $E(f_1)<E(f_2)$. For any $E(f_1)<c<E(f_2), $ since the energy function $E:{\cal M}ap_{k,p}\rightarrow {\bf R}$ is continuous and conformal invariant, hence  $G$-invariant, the inverse images $E^{-1}((-\infty, c))$ and $E^{-1}((c, \infty))$, denoted by  ${W}(f_1)$ and $ { W}(f_2)$,  are two open
  $G$-sets in ${\cal M}ap_{k,p}$ without intersection, containing $f_1$ and $f_2$ respectively. 

  \QED

 ${\bullet}$    Part III: the case that   both  $f_1$   and  $f_2$  are trivial.

   Let $B_{\epsilon'_1}(c_1)$ and  $B_{\epsilon'_2}(c_2)$ be two open balls in $M$, which do not intersect. Here $c_1$ and $c_2$ are the values of the two constant maps $f_1$ and $f_2$ respectively. Clearly if $||h_i-f_i||_{C^0}=max_{x\in \Sigma}|h_i(x)-c_i|<\epsilon'_i, i=1, 2, $ then the image of $h_i$ is contained in $B_{\epsilon'_i}(c_i).$ Moreover, since for  any  $h_i$ and $g_i\in G$,
  the  image of $h_i\circ g_i=$ the  image of $h_i,$ for any  $h_1$ and $h_2$ as above, their $G$-orbits
  $G\cdot h_1$ and  $G\cdot h_2$ do not intersect. Clearly by our assumption   for $\epsilon_i<<
  \epsilon'_i,$
  any  $h_i, i=1, 2$ in $U_{\epsilon_i}(f_i)$ satisfies the condition $||h_i-f_i||_{C^0}<\epsilon'_i$, hence $G\cdot U_{\epsilon_1}(f_1)$
  and $G\cdot U_{\epsilon_2}(f_2)$ do not intersect.

  \QED

\section{ Isotropy groups of the  weakly stable $L_k^p$-maps}
Let  $f$ be a weakly stable $L_k^p$-map.
Consider  a nontrivial unstable component  $f_v:\Sigma_v\simeq {\bf C P }^1\rightarrow M$. The reparametrization group $G_{\Sigma_v}$ of $f_v$ is the group $G_i, i=0, 1, 2$
 considered as subgroup of  ${\bf PSL}(2, {\bf C})$ that fixes the $i$ points that corresponds to the $i$ double points of $\Sigma_v\simeq {\bf C P }^1$.

\begin{lemma}

If $f_v:\Sigma_v\simeq {\bf C P }^1\rightarrow M$  is  not trivial, the identity component $\Gamma^0_{f_v}$ of the isotropy  group $\Gamma_{f_v}$ in $G_i, i=0, 1, 2$ is the standard $S^1$ in $SO(3)\subset G_i\subset {\bf PSL}(2, {\bf C})$ upto a conjugation in $G_i$ if $\Gamma^0_{f_v}$ has a positive dimension.

\end{lemma}

\proof

We only need deal with $G_0$ and $G_1$.

Consider $G_0={\bf PSL}(2, {\bf C})$ first. Since $f_v$ is nontrivial, (i) it is weakly stable so that $\Gamma_{f_v}$ is compact and (ii) the action of $G_0=\Gamma_{f_v}$ on $\Sigma_v$ can not be transitive
so that the connected  maximal compact subgroup $SU(2)/{\bf Z}_2\simeq SO(3)$ of ${\bf PSL}(2, {\bf C})$  and its conjugations are
not contained in $\Gamma_{f_v}$. Since  the groups in the conjugate class of $ U(2)$ are  the  all possible connected maximal compact subgroups of $GL(2, {\bf C})$, so are the groups in the conjugate class of $ SU(2)/{\bf Z}_2$ for ${\bf PSL}(2, {\bf C})$.
Hence   upto a conjugation $\Gamma^0_{f_v}$ is contained in $SO(3)$ as a nontrivial proper compact subgroup. It is well-known that in this case
$\Gamma^0_{f_v}$ is the standard $S^1$ in $SO(3)$ upto a conjugation.

 For  the case that $f_v:\Sigma_v\rightarrow M$ has one double point $\infty \in {\bf P}^1= {\bf C}\cup \{\infty\}$, since  $\Gamma_{f_v}$ is compact, any element $g\in \Gamma_{f_v}\subset G_1$ has the form
 $g(z)
=az+b$ with $a\not =0, 1.$    
Since any element in $\Gamma_{f_v} $ here has two fixed points, given two nontrivial elements $g_1$ and $g_2$  in $\Gamma^0_{f_v}$, we may assume that, $0$ is the other  fixed point for $g_1$ so that $g_1(z)=a_1 z$.
 If $g_2(z)=a_2z+b$ with $b\not =0$,
then $g_3(z)=g_1^{-1}
\circ g_2\circ g_1(z)=a_2z+a_1^{-1}b.$ Hence $g_3\circ g_2^{-1}(z)=a_2(a_2^{-1}(z-b))+a_1^{-1}b=z+a_1^{-1}b-b=z+c$ with $c=a_1^{-1}b-b=(a^{-1}_1-1)b\not =0$ by the assumption that $a_1\not= 1$ and  $b\not = 0.$ Then all the translations by $kc$ are  contained in $\Gamma^0_{f_v}$ which contradicts to the compactness of $\Gamma^0_{f_v}$.
Thus all the elements in  $\Gamma^0_{f_v}$ have the two common fixed points, say $0$ and $\infty$ so that  $\Gamma^0_{f_v}$ is the standard $S^1$ in $G_1$ upto a conjugation.

\QED

\begin{pro}
	If all unstable components of a weakly stable map are
	$2$-dimensional,  it is stable.
\end{pro}
\proof

If a nontrivial unstable component  $f_v:\Sigma_v\simeq {\bf C P }^1\rightarrow M$ of $f$  is   $2$-dimensional at a point $x_0$, then $f_v$ is a local embedding on a small $\epsilon$-disk $D_\epsilon(x_0)$ centered at $x_0$.

By the lemma above,   if the identity component ${\Gamma}^0_{f_v}$  is  nontrivial, it  is $\simeq S^1$ and  can be conjugated into the maximal subgroup of $K=SU(2)/\{\pm 1\}=SO(3)$ of  $G_0$ or the standard $S^1$ in $G_1$ or $G_2$. In term of actions on ${\bf C P }^1\simeq S^2$,   in all three cases, it  is  the standard
$S^1\subset SO(3)$ acting as the  rotations of $ S^2$. 

 By taking a smaller disk in $D_\epsilon(x_0)$,  we may assume that the two fixed points $p_{\pm}$  (the poles) of $S^2$ are not in $D_\epsilon(x_0)$.

Then for $\phi\not={\it id}\in S^1$ but sufficiently close to ${\it id}$, the image  $\phi (D_{\epsilon/2}(x_0))$ is different from $D_{\epsilon/2}(x_0)$ but still inside  $D_{\epsilon}(x_0)$. Since $f_v$ is an embedding on $D_{\epsilon}(x_0)$, the images $f_v(D_{\epsilon/2}(x_0))$ and  $f_v(\phi(D_{\epsilon/2}(x_0)))=f_v\circ \phi(D_{\epsilon/2}(x_0))$ are different. Hence on $D_{\epsilon/2}(x_0)$, $f_v\circ \phi\not=f_v$ already so that $\phi$ is not lying inside   the isotropy group of $f_v$. This is a contradiction.

\QED

\begin{lemma}
	If $f_v: \Sigma_v\simeq {\bf C P }^1\rightarrow M$ is an  unstable  component of a weakly stable map of class $C^1$ with $(\Gamma_{f_v}^T)^0$ being non-trivial, then upto a conjugation by the action of  $G_{\Sigma_v}$,
	$f_v={\bar f}_v\circ \mu$. Here ${\bar f}_v:I=[-1, 1]\rightarrow M$ is a $C^1$-map,  and $\mu:\Sigma_v\simeq  {\bf C P }^1\rightarrow I\subset {\bf R}^1\simeq {\bf so}(2)$  is the "height" function.

\end{lemma}

\proof

By the assumption of the lemma, upto a conjugation $(\Gamma_{f_v}^T)^0\simeq $ the standard $S^1$.
Let $\mu: {\bf C P }^1\rightarrow  {\bf so}(2)$  be the  moment map for the standard action of $S^1$ on ${\bf CP}^1.$
Then upto a conjugation, we may assume that the action of $(\Gamma_{f_v}^T)^0$ has the moment map $\mu$. 

Let $p_{\pm}\in S^2$ be the two poles of $S^2$ that are the two critical points of $\mu.$ Fix a great circle $C\subset S^2$ containing $p_{\pm}$ with the parametrization $[-\pi/2, 3/2\pi)$ such that $p_{-}$ and $p_{+} $ has the coordinates $0$ and $\pi$ respectively.  Let $I_0=(-\epsilon, \pi+\epsilon)$ for a sufficiently small positive $\epsilon$ and $C_0\subset C$  be the corresponding arc. Then the restriction $f_v|_{C_0}:C_0\rightarrow M$ is of class $C^1$. It gives rise  a corresponding $C^1$-smooth map ${\bar f}_v:I=(-1-\delta, 1+\delta)\simeq I_0\simeq C_0 \rightarrow M$. Here the identification $(-1-\delta, 1+\delta)\simeq I_0$ maps
$I=[-1, 1]
$  onto $[0, \pi]$. Then ${\bar f}_v$ so defined satisfies the requirement of the lemma.

\QED

\begin{definition}
A non-trivial $L_k^p$-map $f:S^2\rightarrow M$ is said to be standard $S^1$-invariant (one-dimensional) map if it has the  form in the above lemma. 
\end{definition}

\begin{pro}
	If none of the  unstable  component of a weakly stable map is standard $S^1$-invariant, then it is stable.
	
\end{pro}

Thus among all  weakly stable maps, almost all of them are stable except the special case that there is a  standard $S^1$-invariant unstable component (upto a conjugation).

 We have essentially classified the continuous part of the isotropy group of a weakly stable map. A temptation to classify the discrete part of the isotropy group  was made in the earlier version of this paper \cite{7}.


\begin{thebibliography}{0}
	

	
	\bibitem{3} Hofer, Wysocki and Zehnder: \textsl{Applications of Polyfold Theory I: The Polyfolds of Gromov-Witten Theory,}
	{\it Preprint} (2011), arXiv:math.11072097[math.SG];
	Memoirs of the American Mathematical Society
	2017; 218 pp.
	
	\bibitem{2}
	Kontsevich,	Maxim: \textsl{Enumeration  of rational curves via torus actions}
	{\it Preprint} (1995), arXiv:hep-th.9405035.
	
		
		

	
	
	
	\bibitem{5} Liu, G: \textsl{ Higher-degree smoothness of perturbations I}, Preprint (2016).
	
	\bibitem{6}
	 Liu, G: \textsl{On the Action of Reparametrization Group on the Space of $L_k^p$-maps I},
	{\it  Preprint} (2013). arXiv:1312.4222v3  [math.SG].
	
	\bibitem{7}
	Liu, G: \textsl{On  Stability of  Nodal $L_k^p$-Maps I}, 	
	arXiv: math.1509.07187[math.SG], 2015.	
	
	
	\bibitem{8}
	Liu, G: \textsl{$C^{m_0}$-Smoothness of Evaluation Maps},
	{\it  Preprint} (2013) arXiv:1312.2161v1  [math.SG].
	
	\bibitem{9} McDuff, D and Salamon, D :  \textsl{J-holomorphic Curves and Symplectic topology }, Colloquium Publications, Vol 52, Amer. Math. Soc. Providence, RI, 2012.
	

	
	
	
	
	
	
	
	
\end{thebibliography}
\end{document}